\documentclass[11pt,a4paper]{article}
\usepackage[latin1]{inputenc}
\usepackage[american]{babel}      
\usepackage[T1]{fontenc}
\usepackage{algorithm2e}
\usepackage{amsmath,amssymb,amsthm}
\usepackage{graphicx}
\usepackage{xcolor}   
\usepackage{geometry}
\usepackage{url}
\usepackage{placeins}
\usepackage{listings}
\usepackage{subcaption} 
\usepackage{lipsum} 

\title{{\bf Implementation of neural network operators with applications to remote sensing data}}
         
\author{ {\bf Danilo Costarelli} and {\bf Michele Piconi}\\  
Department of Mathematics and Computer Science \\
            University of Perugia\\
        1, Via Vanvitelli, 06123 Perugia, Italy    \\  
{\small {\tt danilo.costarelli@unipg.it}} - {\tt michele.piconi@unipg.it} }

\date{}

\newcommand{\mau}{\geq}
\newcommand{\miu}{\leq}

\newcommand{\N}{\mathbb{N}}
\newcommand{\R}{\mathbb{R}}
\newcommand{\Z}{\mathbb{Z}}
\newcommand{\disp}{\displaystyle}
\newcommand{\be}{\begin{equation}}
\newcommand{\ee}{\end{equation}}
\newcommand{\phis}{\phi_{\sigma}}
\newcommand{\psis}{\Psi_{\sigma}}

\newcommand{\xx}{\underline{x}}

\newcommand{\uu}{\underline{u}}
\newcommand{\kk}{\underline{k}}

\newtheorem{definition}{Definition}[section]

\newtheorem{theorem}[definition]{Theorem}
\newtheorem{lemma}[definition]{Lemma}

\begin{document}

\maketitle 

\begin{abstract}
In this paper, we provide two algorithms based on the theory of multidimensional neural network (NN) operators activated by hyperbolic tangent sigmoidal functions. Theoretical results are recalled to justify the performance of the here implemented algorithms. Specifically, the first algorithm models multidimensional signals (such as digital images), while the second one addresses the problem of rescaling and enhancement of the considered data. The asymptotic computational complexity of the proposed algorithms is also analyzed. Several applications of the NN-based algorithms for modeling and rescaling/enhancement remote sensing data (represented as images) are discussed, together with numerical experiments conducted on a selection of remote sensing (RS) images from the (open access) RETINA dataset. A comparison with classical interpolation methods, such as bilinear and bicubic interpolation, shows that the proposed algorithms outperform the others, particularly in terms of the Structural Similarity Index (SSIM).

\end{abstract}

\section{Introduction} \label{sec1}
Artificial neural networks (NNs), inspired by the structure and function of the human brain, were first conceptualized in the early 20th century. They are composed of layers of interconnected nodes, each simulating the behavior of a biological neuron, allowing for the representation and processing of complex patterns. Since their introduction, NNs have become a turning point in artificial intelligence and machine learning, with applications ranging from natural language processing and computer vision to more theoretical disciplines such as Approximation Theory and Numerical Analysis (see, e.g., \cite{CS13, CS13b, Costarelli2014, Costarelli2015,CS-2015, An2016,KS19,LTY19,ZH1,ZH2,KOKR1,Costarelli2022,LMP24}).
\\
An important property of NNs lies in their ability to approximate functions, opening the way for many practical applications \cite{An2011a,An2011, COCOKA1,Kad23,CARS24}. Cybenko's pioneer result \cite{CY} showed the approximation capability of single-layer NNs with sigmoidal activation functions, generating an important subsequent scientific interest. Since this work employed non-constructive tools from Functional Analysis, the necessity arose to develop more constructive methods, such as the study of neural network (NN) operators, which have roots in positive linear operators and Approximation Theory.
\\
In this framework, the early contribution by Cardaliaguet and Euvrard \cite{CAEU1} introduced operators with bell-shaped activation functions, establishing the groundwork for further generalizations. These advancements have been crucial for addressing practical challenges in applications, requiring the property of approximating even not necessarily continuous functions, which model the most of real-world signals.

In this paper, we deal with a family of multidimensional NN operators of Kantorovich type (see equation (\ref{defNNkantorovich}) below, for the definition), introduced in \cite{CS-2014} as a generalization of the classic NN operators. Working in a multidimensional setting allows the application of theoretical results related to these operators to image and signal processing, or, more in general, to real-world data.

From a mathematical viewpoint, indeed, it is well-known that digital images can be represented as two-variable step functions with compact support. Over the years, several models have been introduced for applications in the field of image processing, using algorithms such as the Sampling Kantorovich (SK) algorithm (see, e.g., \cite{costarelli2020comparison}) and fuzzy-type algorithms (see \cite{jurio2011construction, CARS24}), leading to applications in various areas, including medicine and engineering. In particular, the SK algorithm proves to be effective in reducing signal noise by rescaling images and acting as a smoothing filter. This is mathematically justified by the regularizing properties of the integral mean appearing in (\ref{defNNkantorovich}), which is essentially a continuous convolution involving the function.

The primary aim of this paper is to introduce an algorithm for image reconstruction based on the theory of NN Kantorovich operators (\ref{defNNkantorovich}). This goal is supported by several theoretical results, recalled in Section \ref{theory}, which represent the starting point for the implementation presented herein. Among these results is a pointwise and uniform convergence theorem in several variables, provided in Theorem \ref{theorem1}. By using NN Kantorovich operators, we are able to approximate not only continuous functions but also those that are not necessarily continuous, such as digital images. This is shown in Theorem \ref{theorem3}, which establishes an $L^p$-convergence theorem in the multidimensional setting. Furthermore, the theory provided for NN operators is not limited to convergence results only but also provides results on the order of approximation. This is achieved through quantitative estimates based on the modulus of smoothness of the function, as shown in Theorem \ref{theorem3}; this result can be useful for understanding the accuracy of the proposed algorithms' performance.

Section \ref{sec_alg} is entirely devoted to the introduction of the tow novel algorithms based on Kantorovich NN operators activated by the hyperbolic tangent sigmoidal function. Specifically, the main tasks achieved by such algorithms are the following:
\begin{itemize}
    \item Algorithm \ref{alg1}: it aims to provide an approximate model of data (such as RS data);
    \item Algorithm \ref{alg2}: it performs rescaling/enhancement of data (again, such as RS data).
\end{itemize}

Section \ref{CC} is instead devoted to the analysis of the computational complexity of the proposed algorithms. Finally, Section \ref{sec_numerical} presents numerical tests applying both Alg. \ref{alg1} (see Subsection \ref{subsection1}) and Alg. \ref{alg2} (see Subsection \ref{subsection2}) to a selection of satellite images from the official open access RETINA dataset, providing a practical application in the field of remote sensing.
The performances of Alg. \ref{alg1} and Alg. \ref{alg2} are evaluated through numerical experiments and comparisons with classical interpolation methods, such as bilinear and bicubic interpolation. Performances are measured using well-known similarity indexes like the Peak Signal-to-Noise Ratio (PSNR) and the Structural Similarity Index (SSIM) (both recalled in Section \ref{sec_methods}).

\section{Approximation results by NN operators}\label{theory}

 Below we recall some basic facts concerning the Kantorovich NN operators, that have been widely studied in recent years in the frame of Approximation Theory.

 Let $\sigma: \R \to \R$ be a given measurable function; $\sigma$ is called a {\em sigmoidal function} if $\lim_{x \to -\infty}\sigma(x)=0$ and $\lim_{x \to +\infty}\sigma(x)=1$, according to the definition given by Cybenko in 1989 (\cite{CY}).
 \vskip0.2cm

   Now, we introduce the following class of functions, that we denote by the symbol ${\cal D}$, i.e., the class containing any non-decreasing sigmoidal function $\sigma$, with $\sigma(1)<1$, satisfying the following assumptions:
\begin{itemize}
\item[$(D 1)$] $\sigma(x)-1/2$ is odd (symmetric);
\item[$(D 2)$] $\sigma$ belongs to $C^2(\R)$ with $\sigma''(x)\leq 0$ for every $x \in [0,+\infty)$;
\item[$(D 3)$] $\sigma(x) \leq C\, |x|^{-\alpha-1},$ as $x \to -\infty$, for some $C>0$ and $\alpha>0$,
\end{itemize}
as introduced in \cite{CS-2014}.

Let now any $\sigma \in {\cal D}$ be fixed. The density function $\phis$ generated by $\sigma$ can be defined by:
\be \label{density-functions}
\phi_{\sigma}(x)\, :=\, \frac{1}{2}[\sigma(x+1)-\sigma(x-1)], \hskip1cm x \in \R.
\ee
Using the above definition in (\ref{density-functions}), we can immediately recall the following operators introduced in \cite{CS-2014}.
\begin{definition}\label{KNN}
Let $\sigma \in {\cal D}$ be fixed and $n \in \N$.  
We define the Kantorovich NN operators, by:
\be \label{NNopS}
(K_n f)(x)\ =\ {\disp \sum_{k=\lceil na \rceil}^{\lfloor nb \rfloor-1} \left[ n \int_{k/n}^{(k+1)/n} f\left(u  \right)\, du \right] \phis\left( n x - k \right) \over \disp \sum_{k=\lceil na \rceil}^{\lfloor nb \rfloor-1} \phis(nx-k)}, \hskip1cm x \in I:=[a,b],
\ee
 where $f:I \to \R$ is a given measurable and bounded function and $\lceil \cdot \rceil$ and $ \lfloor \cdot \rfloor$ denote the "ceiling" and the "integer part" of a given number. Here, $\phis$ is the density function defined in (\ref{density-functions}).
\end{definition}

Below, we mention the multivariate extension of the operators $K_n$. Denoting by $I^d:=[a_1,b_1] \times ... \times [a_d, b_d] \subset \R^d$ a fixed d-dimensional rectangle of $\R^d$, and by the set:
\be
{\cal V}_n\ :=\ \left\{ \kk \in \Z^d:\, \lceil na_i \rceil \leq k_i \leq  \lfloor nb_i \rfloor-1,\ i\, =\, 1,\,...,\, d \right\},
\ee
we can define the following:
\be\label{defNNkantorovich}
     K^d_n(f,\xx)\ :=\ \frac{\disp \sum_{\kk \in {\cal V}_n} \left[ n^d \int_{R^n_{\kk}} f
           \left(\uu\right) \, d\uu \right]
        \Psi_{\sigma}(n \xx - \kk)}
      {\disp \sum_{\underline{k} \in {\cal V}_n}               \Psi_{\sigma}(n\xx - \kk)},  \quad \xx \in I^d, 
\ee
where:
\be
R^n_{\kk} := 
     \left[ \frac{k_1}{n}, \frac{k_1+1}{n} \right] \times  \cdots \times 
   \left[ \frac{k_d}{n}, \frac{k_d+1}{n} \right], 
\ee
are the rectangles in which we will compute the multiple integrals (means) of the considered function of several variables $f:I^d\to \R$, while:
\be
\Psi_{\sigma}(\xx) := \phi_{\sigma}(x_1) \cdot \phi_{\sigma}(x_2) 
           \cdots \phi_{\sigma}(x_d), 
                  \quad \quad    \xx := (x_1, ..., x_d) \in \R^s.
\ee 
is the multivariate (tensor-product) density function defined by $\sigma \in {\cal D}$.

Concerning the above function $\psis$, the following useful properties have been established in \cite{CS-2014}.
\vskip0.2cm

\begin{lemma} \label{lemma1}
(a)  For any $\xx \in \R^d$, we have $\sum_{\kk \in \Z^d} \Psi_{\sigma}(\xx - \kk) = 1$.

\vskip0.2cm

\noindent (b) The following series $\sum_{\kk \in \Z^d} \Psi_{\sigma}(\xx - \kk)$ it turns out uniformly convergent 
on every compact sets of $\R^d$.

\vskip0.2cm

\noindent (c) Let $\|\cdot \|$ the usual maximum norm of $\R^d$, i.e., $\| \xx \| := \max\left\{ |x_i|, \, i = 1, ..., d \right\}$, with $\xx \in 
\R^d$. For every $\gamma > 0$, we get
$$
   \lim_{n \to +\infty} \sum_{\| \xx - \kk \| > \gamma n} \Psi_{\sigma}(\xx - \kk) = 0,
$$
uniformly with $\xx \in \R^d$. In particular, for every $\gamma > 0$
and for every $0 < \nu < \alpha$, 
$$
   \sum_{\| \xx - \kk \| > \gamma n} \Psi_{\sigma}(\xx-\kk) 
      = \mathcal{O}(n^{-\nu}),    \hskip0.5cm  for  \hskip0.5cm   n \to +\infty,
$$
where $\alpha > 0$ is the parameter in assumption $(D 3)$.

\vskip0.2cm
\noindent (d) For every $\xx \in I^d$, we have:
$$
\sum_{\kk \in {\cal V}_n} \psis(n\xx-\kk)\ =\  \prod^d_{i=1} \sum^{\lfloor nb_i \rfloor-1}_{k_i = \lceil na_i \rceil}
       \phi_{\sigma}(nx_i - k_i)\ \mau\ [\phi_{\sigma}(2)]^d\ >\ 0.
$$ 

\vskip0.2cm
\noindent (e) It turns out that $\psis$ belongs to $L^1(\R^d)$, with $\|\psis\|_1=1$, where $\| \cdot\|_1$ denotes the usual $L^1$-norm of a given summable function.
\end{lemma}

We now concentrate the attention to the approximation properties of $K^d_n f$, $d \geq 2$. Hence, using the properties outlined in Lemma \ref{lemma1}, we can recall the following pointwise, uniform, and $L^p$ convergence results (see \cite{CS-2014} again), that will be very useful in the next sections.

\begin{theorem}  \label{theorem1}
  Let $f: I^d \to \R$ be a given bounded function. Hence the operators $K^d_n f$ turns out to be well-defined and, furthermore we also have: 
$$
     \lim_{n \to +\infty} K^d_n(f,\xx) = f(\xx)
$$
at each point $\xx \in I^d$ where $f$ is continuous. Moreover, if $f \in C(I^d)$, i.e., a given continuous function on the whole domain $I^d$, then
\[
    \lim_{n \to +\infty} \sup_{\xx \in I^d} |K^d_n(f, \xx) - f(\xx)|
       = \lim_{n \to +\infty} \|K^d_n(f, \cdot) - f(\cdot) \|_{\infty} = 0.
\]
\end{theorem}

The second part of the above theorem is crucial in order to establish the $L^p$ approximation theorem, with $1 \miu p <+\infty$, that will be recalled below.

First, we can observe that, the operators are well-defined also if $L^p$ functions are taken into account. More precisely:
\begin{theorem} \label{theorem2}
The following inequality: 
$$
        \| K^d_n(f, \cdot)\|_p 
    \miu \frac{1}{[\phi(2)]^{d/ p}} \, \| f \|_p,
$$
is satisfied, for any $f \in L^p(I^d)$, $1 \miu p < +\infty$, where $\| \cdot \|_p$ is the usual $L^p(I^d)$ (integral) norm.
\end{theorem}

On the bases of the proof of the above inequality there is the application of the well-known Jensen inequality for convex functions (see \cite{CS-2015}). Now, we can finally recall the following crucial result.

\begin{theorem} \label{theorem3}
For any given $f \in L^p(I^d)$, $1 \miu p <+\infty$, we get
$$
    \lim_{n \to +\infty} \| K^d_n(f, \cdot) - f \|_p = 0.
$$
\end{theorem}

Other than the above convergence results, also the problem of the order of approximation has been considered for the above family of positive linear operators. For a comprehensive study we refer to the following recent papers \cite{COCOKA1,COCO1}, in which both the cases of the uniform and $L^p$-norm have been considered. In this sense, below we mention the following very complete result occurring with respect to the uniform norm. The same result can be easily formulated for functions in $L^p$ (see \cite{CP1} for the one-dimensional theory).

\begin{theorem} \label{th4}
Let $\sigma \in {\cal D}$ and $f\in C\left( I^d\right)$ be fixed. Denoting by:
$$
\omega (f,\delta )\, :=\, \sup \left\{ \left\vert f(\mathbf{\underline{x}})-f(%
\mathbf{\underline{y}})\right\vert :\, \mathbf{\underline{x}},\, \mathbf{%
\underline{y}} \in I^d,\, \left\Vert \mathbf{\underline{x}}-\mathbf{\underline{%
y}}\right\Vert _{2}\leq \delta \right\} \text{,} \quad \delta>0,
$$
the usual modulus of continuity of a given function, where $\|\cdot\|_2$ denotes the usual Euclidean norm of $\mathbb{R}^d$, we can state what follows.

\vskip0.2cm

\noindent (I)  If $\sigma$ satisfies condition $(D 3)$ for $0<\alpha < 1$, then there exists a constant 
$C>0$, such that
\begin{equation*}
\left\vert K^d_{n}(f,\mathbf{\underline{x}})-f(\mathbf{\underline{x}}%
)\right\vert\, \leq\, C\omega \left( f,\,\frac{1}{n^{\alpha}}\right) \text{,
for all }{\underline{x}}\in I^d,
\end{equation*}%
and $n\in \mathbb{N}$.
\vskip0.2cm

\noindent (II)  If $\sigma$ satisfies condition $(D 3)$ for $\alpha > 1$, then there exists a constant 
$C>0$, such that
\begin{equation*}
\left\vert K^d_{n}(f,\mathbf{\underline{x}})-f(\mathbf{\underline{x}}%
)\right\vert\, \leq\, C\omega \left( f,\,\frac{1}{n}\right) \text{,
for all }{\underline{x}}\in I^d,
\end{equation*}%
and $n\in \mathbb{N}$.

\vskip0.2cm

\noindent (III)  If $\sigma$ satisfies condition $(D 3)$ for $\alpha = 1$, then there exists a constant 
$C>0$, such that
\begin{equation*}
\left\vert K^d_{n}(f,\mathbf{\underline{x}})-f(\mathbf{\underline{x}}%
)\right\vert\, \leq\, C\omega \left( f,\,\frac{\ln n}{n}\right) \text{,
for all }{\underline{x}}\in I^d,
\end{equation*}%
and $n\in \mathbb{N}$.
\end{theorem}

Finally, we remark that the theory of neural network operators represents a very active field; in this sense we can also mention the following very recent extension of the above definition, that currently is available only in the one-dimensional case:
\begin{equation}
\left( D_{n}^{\sigma ,\chi }f\right) (x):=\frac{\displaystyle\sum_{k=\lceil {na \rceil }}^{\lfloor {nb \rfloor }-1} 
\left[n\, \int_{a}^{b}\chi (nt-k)f(t)\, dt\right] \phi _{\sigma }(nx-k)}{\displaystyle\sum_{k=\lceil {na \rceil }}^{\lfloor {nb \rfloor }-1}
\left[n\, \int_{a}^{b}\chi (nt-k)\, dt\right] \phi _{\sigma }(nx-k)}, \quad x\in \lbrack a,b].
\end{equation}
for $\sigma \in {\cal D}$ and where $\chi :\mathbb{R}\rightarrow \lbrack 0,+\infty )$ be bounded and $L^{1}$%
-integrable on $\mathbb{R}$, such that
\begin{equation}
\int_{0}^{1}\chi (u)\,du\ =:\ {\cal K}\ >0,
\label{unity}
\end{equation}
and its \textit{discrete absolute moment of order $0$} is finite, i.e.,
\begin{equation*}
M_{0}(\chi):=\sup_{u\in\R} \sum_{k\in\Z} \chi(u-k)<+\infty.
\end{equation*}
The operators $D_{n}^{\sigma ,\chi }$ are called the Durrmeyer-type neural network (NN) operators associated to $f$ with respect to $\phi_{\sigma }$ and $\chi$ and have been introduced in \cite{CCNP1}. Very recently, also a Steklov-type generalization of NN operators (inspired to \cite{CO8}) has been considered in \cite{Acar-2024}. Extensions of the above operators to the multidimensional case should be desirable, but their are currently not available.


\section{Implementation of the operators, pseudo-code of the proposed algorithms, and considered sigmoidal functions}\label{sec_alg}

In this section we provide a new application of the above Kantorovich NN operators to the setting of data (image) modeling, data rescaling and enhancement. The proposed algorithms will be used to study applications to RS data, according with the research project named "RETINA", mentioned in the funding section of the present paper.

In order to achieve the desired result, we will consider the operators $K^d_n$ for $d=2$. First, we recall that any gray scale image (i.e., a given matrix of pixels $ A=(a_{i,j})_{\substack{i=1,\dots,M \\ j=1,\dots,N}}$ of dimension $M \times N$) can be modeled as a function with support contained in the bi-dimensional rectangle $[0,M] \times [0, N] \subset \R^2$, such that in any open square $(i-1,i)\times(j-1,j)$ of side-length equal to one, it assume a constant value coinciding with the pixel luminance $a_{i,j}$, $i=1,...,M$ and $j=1,...,N$. In what follows, we refer to such kind of functions as {\em image function} associated to $A$, formally defined by:
\be \label{calA}
{\cal A}(\xx)\, =\, {\cal A}(x,y)\ =\ \sum_{i=1}^M \sum_{j=1}^N a_{i,j}\, \chi_{i,j}(\xx), \quad \xx=(x,y) \in I^2_{\cal A}:=[0,M]\times[0,N],
\ee
where $\chi_{i,j}(\xx)$ denote the characteristic functions of the sets $(i-1,i]\times(j-1,j]$, $i=1,...,M$ and $j=1,...,N$.

From the above definition of the function ${\cal A}$ we can deduce that pointwise reconstructions, as well as, $L^p$ convergence results can be deduced for the operators $K^2_n {\cal A}$ to ${\cal A}$ itself, as $n \to +\infty$.

Hence, by the Kantorovich NN operators $K^2_n$ applied to the function ${\cal A}$, we are able to provide (approximate) model of data (images), and to perform rescaling/enhancement in term of their resolution. Below, we provide the pseudo-codes of the proposed algorithms (see Alg. \ref{alg1} and Alg. \ref{alg2}). We begin with Alg. \ref{alg1} for data modeling.

\RestyleAlgo{ruled}
\begin{algorithm}[hbt!]
\caption{Pseudocode of the algorithm for data modeling based on Kantorovich NN operators}\label{alg1}
\textbf{Goal:} Modeling of the original data (image) ${\cal A}$ using the NN operators activated by $\sigma \in {\cal D}$.\\
\textbf{Input data:} Original data ${\cal A}$ ($M\times N$ dimension); parameter $n\in\mathbb{N}$.\\
\begin{description}
\item[-] Definition of the multivariate function $\Psi_{\sigma}$;
\item[-] Construction of matrices ${\cal M}$ of the necessary mean values (multiple-integrals);
\item[-] Computations of the vector-arguments of the multivariate density function $\Psi_{\sigma}$;
\item[-] Implementation of a grid of nodes on the set $I^2_{\cal A}$ useful to establish the modeling.
\item[-] Iteration:
\end{description}
 \For{$i=1,\dots,M$ and $j=1,\dots,N$}{
  sum over $\kk \in {\cal V}_n$ of all meaningful terms of the form $\Psi_\sigma(n{\tt{x}}-{\tt{k}})$ convoluted with the matrix of the mean values, where $\xx$ is assumed on the constructed grid of nodes\;
  }
\vspace{0.3cm}
\KwResult{The modeled data by the Kantorovich NN operators.}
\end{algorithm}

At the same time, the above operators can be used for image rescaling and enhancement. In this case the pseudo-code of the algorithm can be found in Alg. \ref{alg2}.

\RestyleAlgo{ruled}
\begin{algorithm}[hbt!]
\caption{Pseudocode of the algorithm for data rescaling and enhancement based on Kantorovich NN operators}\label{alg2}
\textbf{Goal:} Rescaling and enhancement of the original data (image) ${\cal A}$ using the NN operators activated by $\sigma \in {\cal D}$ with a certain (final) resolution based on the scaling factor S.\\
\textbf{Input data:} Original data ${\cal A}$ ($M\times N$ dimension); parameter $n\in\mathbb{N}$.\\
\begin{description}
\item[-] Definition of the multivariate function $\Psi_{\sigma}$;
\item[-] Setting of the final size of the reconstructed image: $(M \cdot S) \times (N \cdot S)$;
\item[-] Construction of matrices ${\cal M}$ of the necessary mean values (multiple-integrals);
\item[-] Computations of the vector-arguments of the multivariate density function $\Psi_{\sigma}$, scaled by the parameter $S$;
\item[-] Implementation of a grid of nodes on the set $I^2_{\cal A}$ useful to establish the scaled/enhanced data.
\item[-] Iteration:
\end{description}
 \For{$i=1,\dots,M$ and $j=1,\dots,N$}{
  sum over $\kk \in {\cal V}_n$ of all meaningful terms of the form $\Psi_\sigma(n{\tt{x}}-{\tt{k}})$ convoluted with the matrix of the mean values, where $\xx$ is assumed on the constructed grid of nodes\;
  }
\vspace{0.3cm}
\KwResult{The scaled/enhanced data by the Kantorovich NN operators.}
\end{algorithm}

Now, a crucial aspect in the application of the above introduced algorithms is represented by the choice of the sigmoidal function $\sigma \in {\cal D}$. In the literature, we can found several instances of functions $\sigma$ (see \cite{CC09}), such as the well-known logistic function (see, e.g., \cite{COCOKA1}) or the ramp function (see, \cite{Costarelli2014}).

Through this paper we mainly deal with the sigmoidal function generated by the hyperbolic tangent function, namely:
$$
\sigma_h(x) := (\tanh x + 1)/2, \quad x \in \R,
$$
where we recall that:
$$
\tanh x\, :=\, {e^x - e^{-x} \over e^x + e^{-x}}, \quad x \in \R.
$$
In particular, other than $(D 1)$ and $(D 2)$, $\sigma_{h}$ satisfies condition $(D 3)$ for every $\alpha > 0$, in view of its 
exponential decay to zero as $x \to -\infty$, hence the corresponding order of (pointwise and uniform) approximation can be deduced from (II) of Theorem \ref{th4}. At the same time, a similar order of approximation with respect of $L^p$-norm can be deduced.


\section{Analysis of the computational complexity of the proposed algorithms} \label{CC}

Since the computational complexity (CC) of both the above algorithms are substantially the same, we consider in details only Alg. \ref{alg2}. Using the same symbols and notations considered in the above pseudo-code we proceed by a step-by-step asymptotic estimate of the computational complexity. In details:
\begin{enumerate}
    \item Input of the original data: the CC is ${\cal O}(M N)$;
    \item Construction of ${\cal M}$: the CC is ${\cal O}(M N n^2)$, where $n$ denotes the parameter of the considered NN operator $K^2_n$;
    \item Computations of the vector-arguments of $\Psi_{\sigma}$ and of the grid of nodes: the CC is ${\cal O}(M N n^2/S^2)$, where $S$ is the scaling parameter;
    \item Number of iteration: the CC is ${\cal O}(M N n^2)$;
    \item CC of each iteration: ${\cal O}(n^2/S^2)$.
\end{enumerate}

In coclusion, based on the above analysis, we can conclude that the full CC of the algorithm is the following:
\begin{equation}
{\cal O}(N M n^4 /S^2).
\end{equation}
It is clear that the estimated CC is quite large; however,as will be showed in the next sections, this is justified by a high accuracy.


\section{Methods: other algorithms for numerical comparison and similarity indexes for the evaluation of the results}\label{sec_methods}

In the next section, we provide some numerical test concerning the performances of the above mentioned Alg.s \ref{alg1} and \ref{alg2} activated by the sigmoidal function $\sigma_h$ generated by the hyperbolic tangent.
\\
In order to evaluate the performances of the proposed algorithms, we will provide a comparison using some very classical interpolation methods of digital image processing, namely the bilinear and the bicubic interpolation. These methods are both well-established and readily accessible, with implementations available in various software and dedicated commands supported by most commonly used programming languages. In particular, to evaluate the performances of Alg.s 1 and 2 considered in the numerical tests, we use the images taken from the official RETINA dataset, described in the next section. 
\\
To assess the results of the numerical tests, we will rely on two widely used similarity indices found in the literature. Denoting by $A$ and $B$ two images with dimensions $N, M \in \mathbb{N}$, the first index we recall is the well-knwon \textit{Peak Signal-to-Noise Ratio} (PSNR), which is defined via the MSE as follows:
\begin{equation}\label{psnr}
\text{PSNR} := 20 \cdot \log_{10} \left( \frac{\max{I}}{\sqrt{\text{MSE}}} \right),
\end{equation}
where \( \max{I} \) is the maximum signal measure in the original image (for 8-bit grayscale images, \( \max{I} = 255 \)), and MSE denotes the \textit{mean square error}, a commonly used metric, defined as:
\[
\text{MSE} := \frac{1}{N \cdot M} \sum_{i=1}^{N} \sum_{j=1}^{M} |a_{ij} - b_{ij}|^2.
\]
It is well known that MSE quantifies the average squared difference between corresponding pixel values of two images $A$ and $B$. Hence, we can say that he PSNR value increases with greater "similarity" to the original image.
\\ 
Another metric is represented by the well-known \textit{Structural Similarity Index} (SSIM) \cite{SSIM2}. Unlike MSE-based measures, which compares pair's differences pixel-by-pixel, this metric performs a comparison between a reference image $A$ and a potentially corrupted version of the same image $B$, based on three independent components extracted at a single spatial scale (resolution): luminance, contrast and structure.
In particular, the luminance comparison is given by
\begin{equation*}
  l(A, B) = \frac{2\mu_A\mu_B + C_1}{\mu_A^2 + \mu_B^2 + C_1}, 
\end{equation*}
where $\mu$ represents the luminance information and $C_1$ is a stabilizing constant, usually set as $C_1=(0.01\times L)^2$ (where $L = 255$ denotes the dynamic range of pixel values for 8-bit images). Meanwhile, the contrast is represented trough the use of standard deviation $\sigma$. Hence, the comparison is
\begin{equation*}
    c(A,B)=\frac{2\sigma_A\sigma_B + C_2}{\sigma_A^2 + \sigma_B^2 + C_2},
\end{equation*}
where $C_2 = (0.03 \times L)^2$. \\Furthermore, the structure element is obtained by  using the sample variances and covariance, denoted by $\sigma^2$ and $\sigma_{AB}$, respectively, as follows
\begin{equation*}
    s(A,B)=\frac{\sigma_{AB}+C_3}{\sigma_A\sigma_B+C_3},
\end{equation*}
with $C_3=C_2/2$. The covariance is usually computed as $$\sigma_{AB}=\frac{1}{N\cdot M-1}\sum_{i=1}^N\sum_{j=1}^M(a_{ij}-\mu_A)\cdot(b_{ij}-\mu_B).$$
Finally, the three components are combined into a unique expression that is
\begin{equation}\label{ssim}
\text{SSIM}(A, B) := l(A, B) \cdot c(A, B) \cdot s(A, B) = \frac{2\mu_A\mu_B + C_1}{\mu_A^2 + \mu_B^2 + C_1} \cdot \frac{2\sigma_{AB} + C_2}{\sigma_A^2 + \sigma_B^2 + C_2}.
\end{equation}
SSIM is a decimal value between -1 and 1, where 1 indicates perfect similarity, 0 indicates no similarity, and -1 indicates perfect anti-correlation.
Over the years, its mathematical properties have been extensively analyzed, yielding various theoretical results (see, e.g., \cite{brunet2012geodesics, brunet2011mathematical}). It is also used in applications related to remote sensing (see, e.g., \cite{SSIM_rm, SSIM_rm1}).

\section{Numerical tests for the application to remote sensing data modeling and rescaling/enhancement }\label{sec_numerical}

\begin{figure}[h!]
    \centering
    \includegraphics[width=0.5\textwidth]{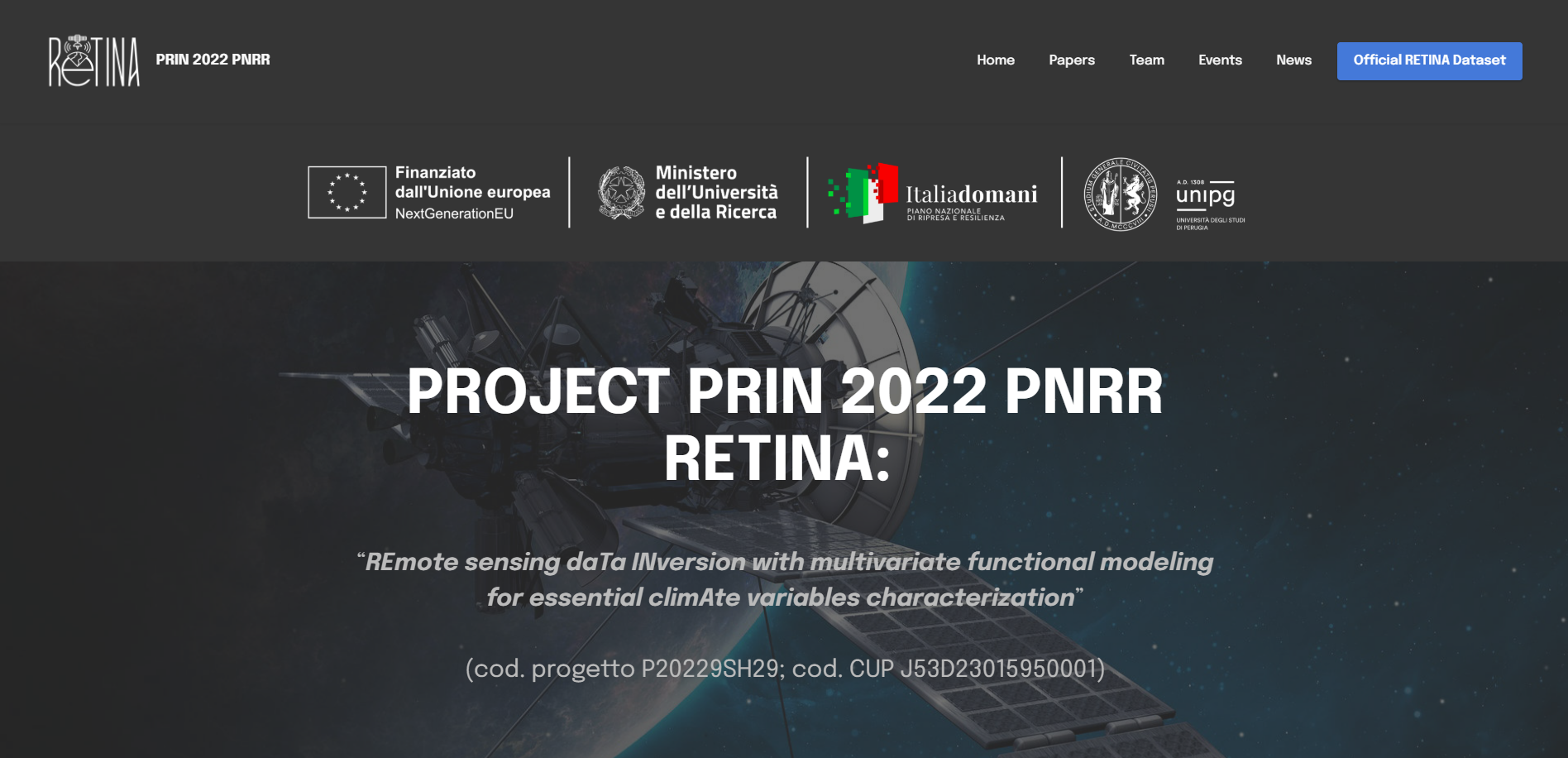} 
    \caption{Screenshot of the RETINA website.}
    \label{cattura_homepage}
\end{figure}

\begin{figure}[h!]
    \centering
    \includegraphics[width=0.5\textwidth]{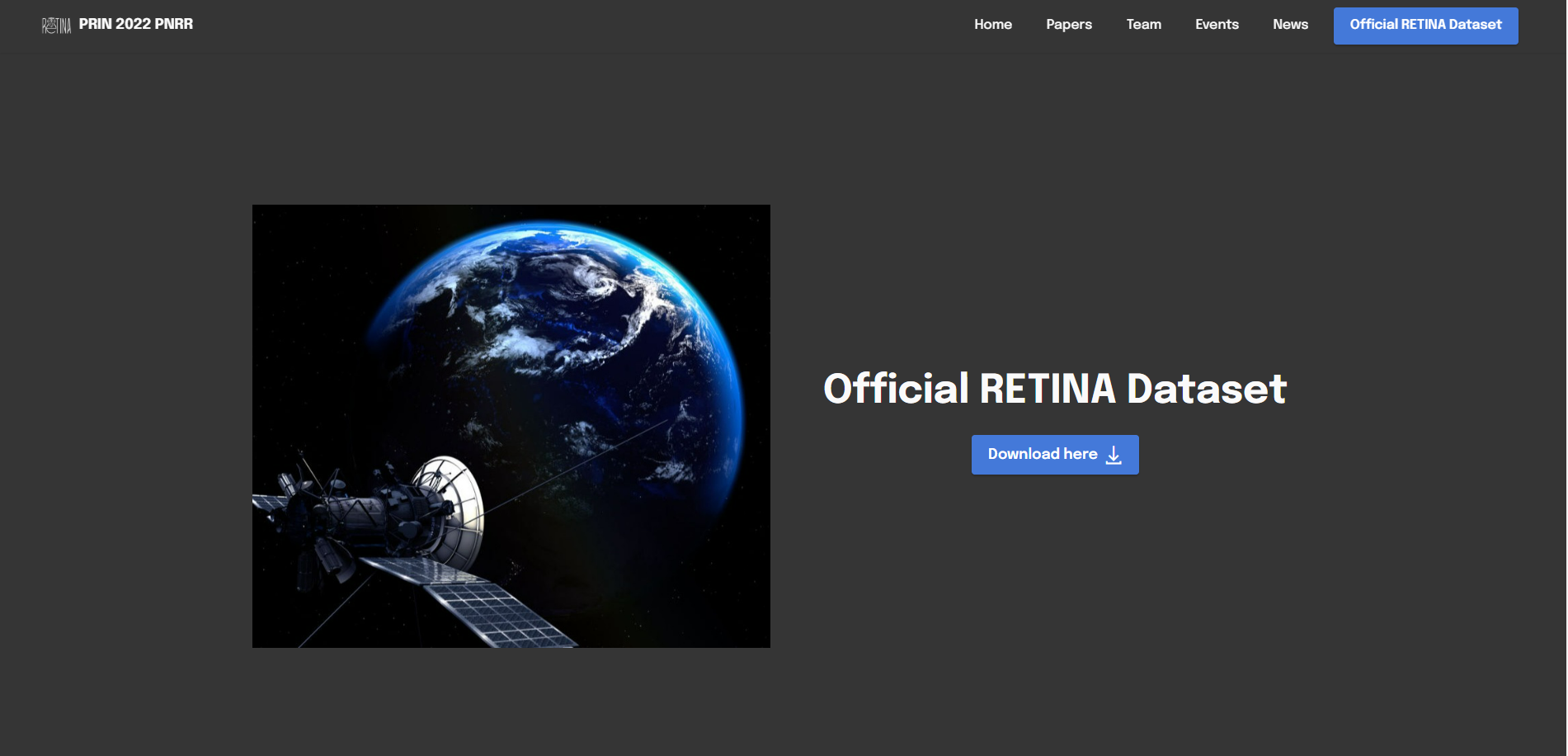} 
    \caption{Screenshot of the dataset from which the images can be freely downloaded.}
    \label{cattura_dataset}
\end{figure}

The motivation for this work arises from the guidelines outlined within the RETINA project (REmote sensing daTa INversion with multivariate functional modeling for essential climAte variables characterization), funded by the Italian Ministry of University and Research (MUR) under the PNRR 2022 call (Figure \ref{cattura_homepage} shows a screenshot of the website \url{https://retina.sites.dmi.unipg.it/index.html}). The RETINA project addresses environmental challenges by developing innovative methods for analyzing data generated from the interaction between electromagnetic waves and the Earth's surface. Specifically, RETINA focuses on key Earth surface variables, such as surface soil moisture (SM) and freeze/thaw (FT) state, which are closely related to Essential Climate Variables (ECVs) used to monitor and understand climate change.
\\
In line with these objectives, the numerical experiments presented in this section aim to validate and evaluate the effectiveness of novel algorithms for image processing based on multidimensional NN operators, which play a fundamental role for the construction of a model for the measured RS data. In particular, the aforementioned algorithms are instrumental in enhancing the understanding of key geophysical variables, thereby supporting RETINA's main goal of connecting theoretical advancements with practical applications in environmental monitoring.
\\ 
The algorithms will be tested on a selection of satellite images collected and made available on the RETINA project website at the following link: \url{https://retina.sites.dmi.unipg.it/dataset.html} (see, Figure \ref{cattura_dataset}). 
\\
All the images available for download on this page are sourced from the Sentinel-1 satellite constellation, which is part of the Copernicus Earth observation component of the European Union Space Programme. The processed images are stored in GeoTIFF format with two bands: VV and VH. These files are compatible with any GIS software, including SNAP and QGIS. Additionally, Python and MATLAB code for reading the images is provided to facilitate their use in the numerical experiments.

To visualize the images in MATLAB, the following commands can be used:

\begin{verbatim}
% Read and process GeoTIFF images in MATLAB
filename_read = 'filenameread.tif';
filename_write = 'filenamewrite.tif';

% Read the GeoTIFF file
[A,R] = readgeoraster(filename_read,'OutputType','double');
info = georasterinfo(filename_read);

% Process the first band (VV)
im1 = A(:,:,1);
im1db = 10*log10(im1);
max_im1db = max(im1db(:));
min_im1db = min(im1db(:));

% Process the second band (VH)
im2 = A(:,:,2);
im2db = 10*log10(im2);
max_im2db = max(im2db(:));
min_im2db = min(im2db(:));

% Display the first band
figure, imagesc(im1db)
colormap('gray'), colorbar, axis image
clim([max_im1db-60, max_im1db])

% Save the image as a GeoTIFF
print(filename_write,'-dtiff','-r600')
\end{verbatim}

\begin{figure*}[h!]
    \centering
    \begin{subfigure}{0.38\textwidth}
        \centering
        \includegraphics[width=\textwidth]{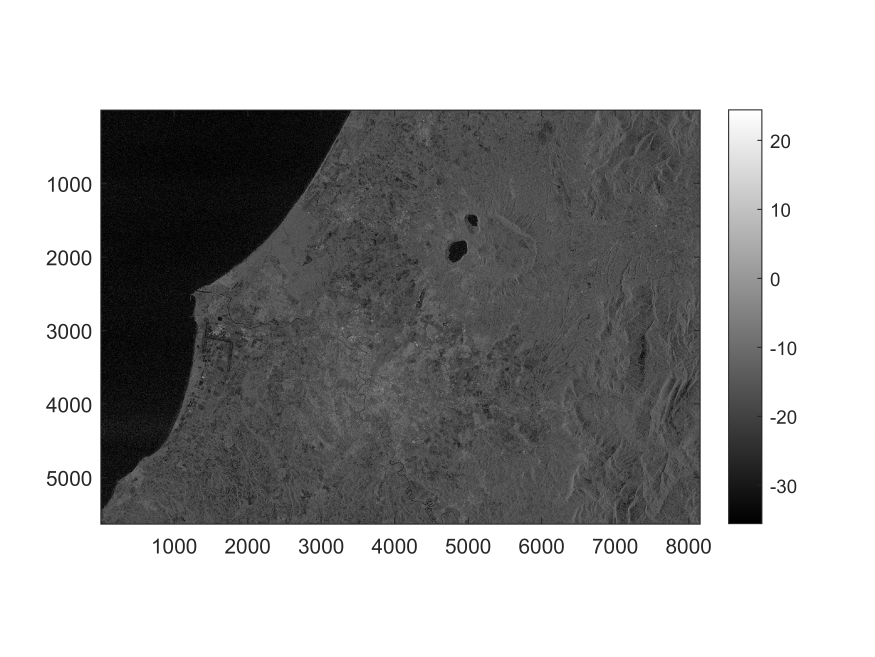} 
        \caption{}
        \label{fig:rome}
    \end{subfigure}%
    \hspace{0.02\textwidth}
    \begin{subfigure}{0.38\textwidth}
        \centering
        \includegraphics[width=\textwidth]{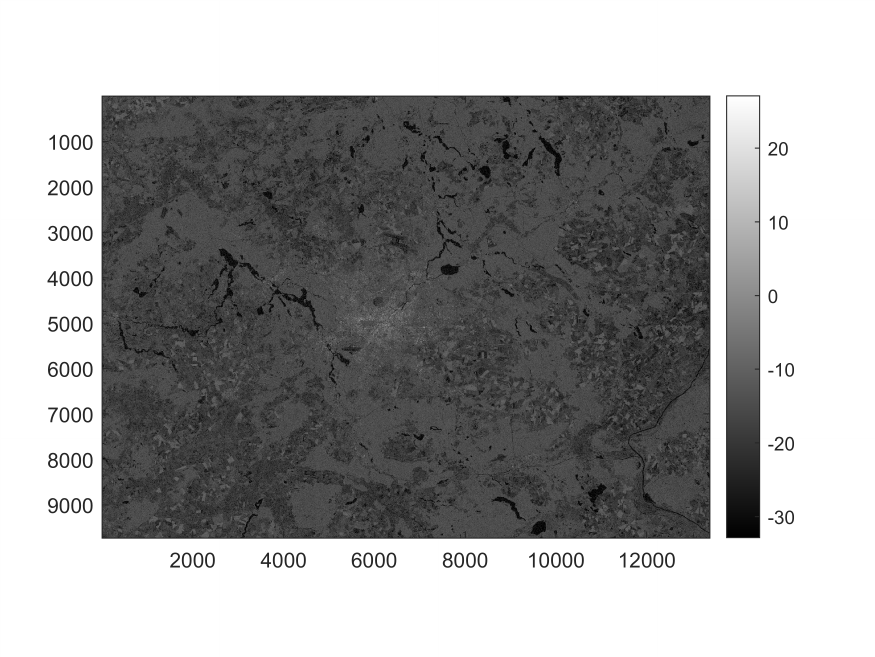} 
        \caption{}
        \label{fig:berlin}
    \end{subfigure}
    
    \vspace{0cm} 
    
    \begin{subfigure}{0.38\textwidth}
        \centering
        \includegraphics[width=\textwidth]{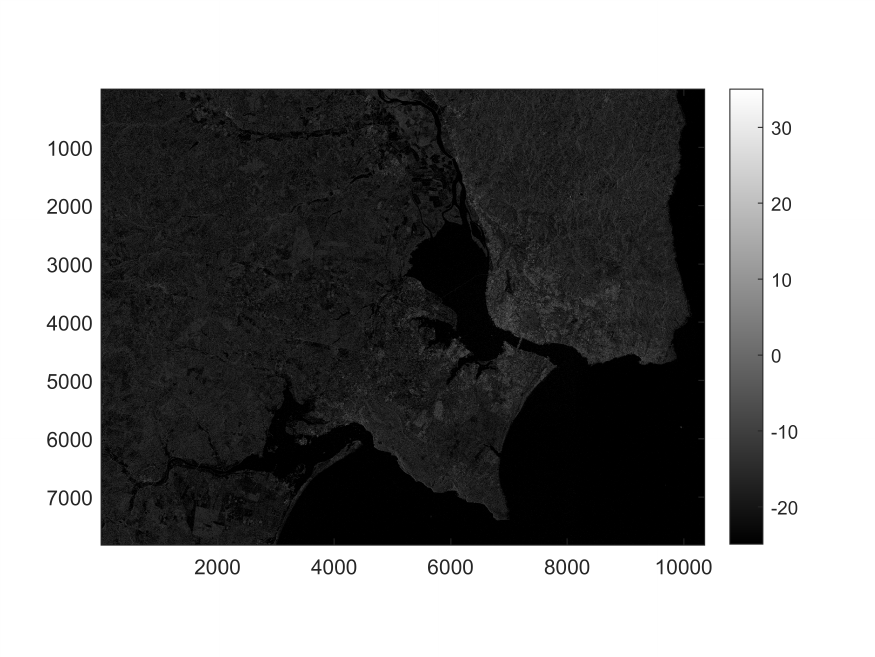} 
        \caption{}
        \label{fig:lisbon}
    \end{subfigure}%
    \hspace{0.02\textwidth}
    \begin{subfigure}{0.38\textwidth}
        \centering
        \includegraphics[width=\textwidth]{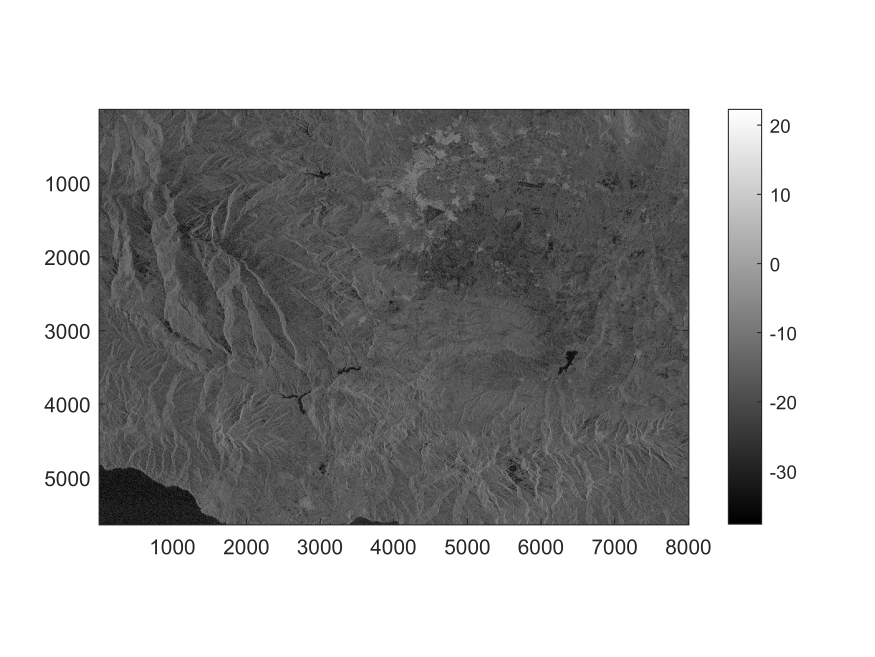} 
        \caption{}
        \label{granada}
    \end{subfigure}
    
    \caption{Satellite images from the RETINA dataset. (a) Rome.tif, (b) Berlin.tif, (c) Lisbon.tif, (d) Granada.tif.}
    \label{satellite_images}
\end{figure*}

To ensure the execution of the MATLAB commands provided, it is necessary to have the \textit{Mapping Toolbox} installed.
\\
We selected four satellite images from the dataset, applying the MATLAB procedure described above. These images will be used to test the NN Alg. 1 and Alg. 2 provided in the previous section. The extracted images are \textit{Rome.tif}, \textit{Berlin.tif}, \textit{Granada.tif} and \textit{Lisbon.tif} and are displayed in Figure \ref{satellite_images}.

\subsection{Applications of Algorithm 1: modeling of RS data}\label{subsection1}

\begin{figure*}[tbph]
    \centering
    \begin{subfigure}{0.22\textwidth} 
        \centering
        \includegraphics[height=3cm]{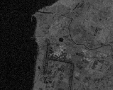}
        \caption{}
        \label{rome2}
    \end{subfigure}
    \hspace{\fill} 
    \begin{subfigure}{0.24\textwidth} 
        \centering
        \includegraphics[height=3cm]{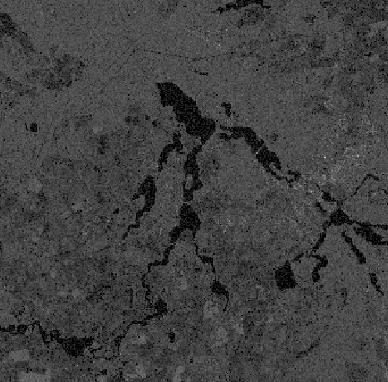}
        \caption{}
        \label{fig:berlin2}
    \end{subfigure}
    \hspace{\fill} 
    \begin{subfigure}{0.27\textwidth} 
        \centering
        \includegraphics[height=3cm]{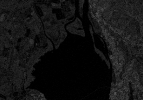}
        \caption{}
        \label{fig:lisbon2}
    \end{subfigure}
    \hspace{\fill} 
    \begin{subfigure}{0.22\textwidth} 
        \centering
        \includegraphics[height=3cm]{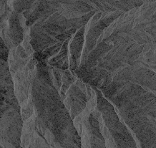}
        \caption{}
        \label{fig:granada2}
    \end{subfigure}
    \caption{Reference images: (a) \textit{Rome}; (b) \textit{Berlin}; (c) \textit{Lisbon}; (d) \textit{Granada}.}
    \label{section_reference_images}
\end{figure*}

Here, we present a series of numerical experiments conducted on selected sections of the satellite images shown in Figure \ref{satellite_images}, which play the role of reference images. The specific sections chosen for each image are displayed in Figure \ref{section_reference_images} and we suppose that their dimension is $M\times N$.
\\
In this subsection, we test Algorithm \ref{alg1} on the selected section of the \textit{Rome} image, shown in Figure \ref{section_reference_images} (a). The goal is to model this image using the NN algorithm with a hyperbolic tangent activation function. 
In particular, Algorithm \ref{alg1}, which focuses on modeling satellite images using neural network operators, is well-suited for the inversion problem. In recent years, real-world applications that involve image analysis have used modeling techniques based on families of quasi-interpolation operators or sampling operators \cite{costarelli2020comparison}. The idea is to represent images (or more generally, data matrices) using functions with compact support, which can then be used to construct the corresponding mathematical operators. This approach has the advantage of providing a model with a known analytical expression, allowing direct control over the parameters. For example, to improve image resolution, we can increase the sampling rate of the operators to achieve the desired result (in the case of NN operators this parameter is represented by $n$). Theoretical analysis provided in Section \ref{theory} justifies these methods. \\ Additionally, since RETINA focuses on retrieving geophysical variables (such as SM and FT), it relies on a modeling framework based on linear operators, such as neural network operators. This framework enables the application of functional analysis techniques to address the inversion problem (see, e.g., \cite{howlett09}).
\\
To evaluate the performance of Alg. \ref{alg1}, we use two similarity indexes: PSNR (see (\ref{psnr})) and SSIM (see (\ref{ssim})).
We outline that, in order to evaluate these indexes in MATLAB, we use its built-in functions \texttt{psnr} and \texttt{ssim}, using the original image in Figure \ref{section_reference_images} (a) as reference. Before performing the calculations, it is necessary to convert the image data from MATLAB's \texttt{uint8} format to \texttt{double}. This step is crucial because applying the \texttt{psnr} function directly to \texttt{uint8} data may produce a zero difference between the original and reconstructed images whenever the difference falls below zero, leading to inaccurate results. This ensures the data type conversion is handled appropriately, preventing inaccuracies in the computed PSNR and SSIM values.
\\
As observed in Table \ref{modelingRometab}, the NN algorithm demonstrates satisfactory performance in terms of quality, as evidenced by both the PSNR and SSIM metrics, even with low values of the parameter $n$. Indeed, for a small value of $n$ ($n=5$), one has a high degree of similarity. 
\vspace{0.3cm}
\begin{table}[h!]
\centering
\renewcommand{\arraystretch}{1} 
\begin{tabular}{|l|c|c|} \hline
 \textbf{$n$-values} & \textbf{PSNR} & \textbf{SSIM} \\ \hline
\textbf{$n=5$} & 43.3192 & 0.9978 \\ \hline
\end{tabular}
\caption{Evaluation of modeling by Alg. \ref{alg1} for the image \textit{Rome} (Figure \ref{section_reference_images} (a)).}
\label{modelingRometab}
\end{table}

Thus, this initial experiment highlights that Alg. \ref{alg1} represents an effective tool for modeling RS data.

\subsection{Applications of Algorithm 2: rescaling/enhancement of RS data}\label{subsection2}
\begin{figure*}[h]
    \centering
    \begin{subfigure}{0.48\textwidth}
        \centering
        \includegraphics[width=\textwidth]{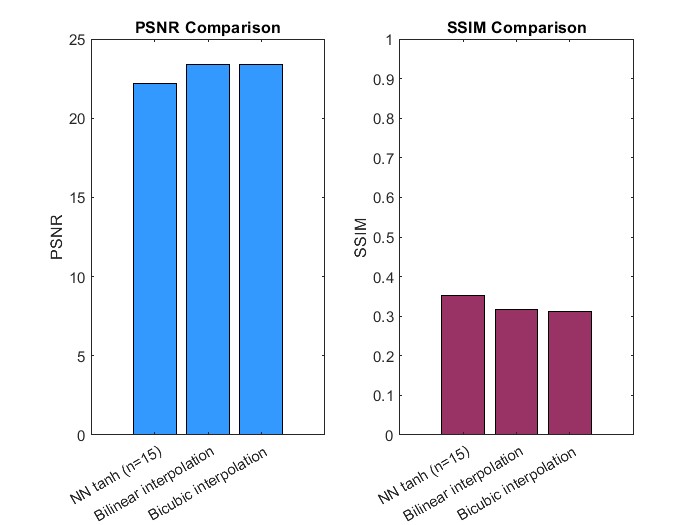} 
        \caption{}
        \label{istorome}
    \end{subfigure}%
    \hspace{0.02\textwidth}
    \begin{subfigure}{0.48\textwidth}
        \centering
        \includegraphics[width=\textwidth]{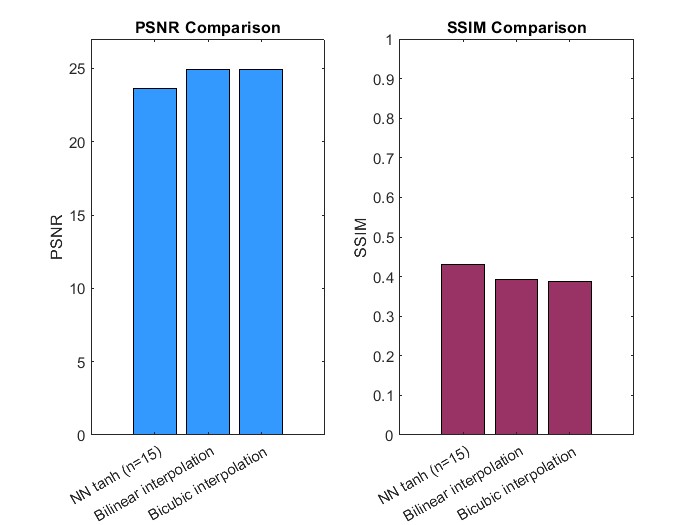} 
        \caption{}
        \label{fig:berlin}
    \end{subfigure}
    
    \vspace{0cm} 
    
    \begin{subfigure}{0.48\textwidth}
        \centering
        \includegraphics[width=\textwidth]{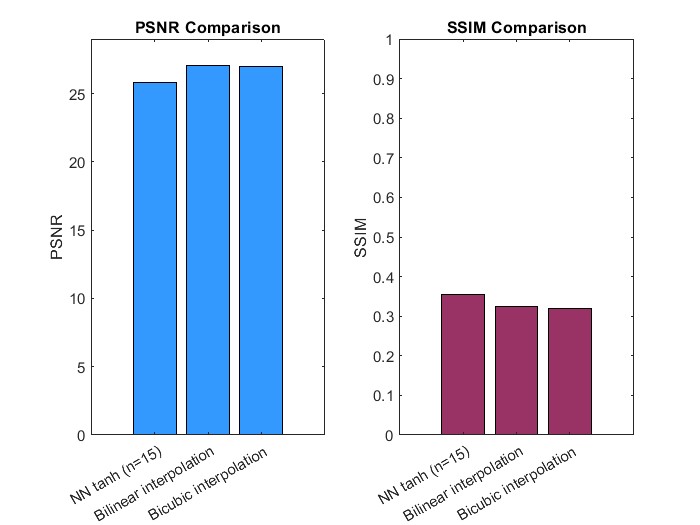} 
        \caption{}
        \label{fig:lisbon}
    \end{subfigure}%
    \hspace{0.02\textwidth}
    \begin{subfigure}{0.48\textwidth}
        \centering
        \includegraphics[width=\textwidth]{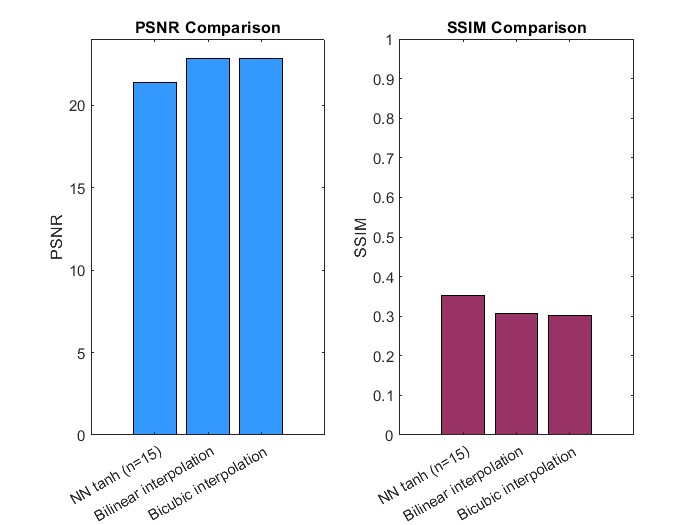} 
        \caption{}
        \label{granada}
    \end{subfigure}
    
    \caption{Comparison of rescaling methods in terms of PSNR and SSIM. (a) \textit{Rome}, (b) \textit{Berlin}, (c) \textit{Lisbon}, (d) \textit{Granada}.}
    \label{istogrammi}
\end{figure*}

Now, we aim to exploit the properties of NN operators activated by hyperbolic tangent function to rescaling/enhancing the selected satellite images. The first step is to downscale them to dimensions $\frac{M}{2} \times \frac{N}{2}$, using the nearest neighbor method (without interpolation) \cite{BD2015}. Following this, the downscaled images were upscaled back to their original size using the methods described above. For the algorithm employing NN operators, namely Alg. \ref{alg2}, we tested values of $n = 5, 10, 15, 20, 25,$ and $30$, where smaller $n$-values correspond to shorter execution times, that in any case, is not significantly high.\\ The numerical results are presented in Tables \ref{Rometab}, \ref{Berlintab}, \ref{Lisbontab}, and \ref{Granadatab}, with the best PSNR and SSIM values highlighted in bold.
\vspace{0.3cm}
\begin{table}[h!]
\centering
\renewcommand{\arraystretch}{1} 
\begin{tabular}{|l|c|c|} \hline
\textbf{Method} & \textbf{PSNR} & \textbf{SSIM} \\ \hline
NN tanh ($n=15$) & 22.2198 & \textbf{0.3522} \\ \hline
Bilinear Interpolation & \textbf{23.4122} & 0.3158 \\ \hline
Bicubic Interpolation & 23.3965 & 0.3118 \\ \hline
\end{tabular}
\caption{Comparison of rescaling methods for the image \textit{Rome} (Figure \ref{section_reference_images} (a)).}
\label{Rometab}
\end{table}
\vspace{0.3cm}
\begin{table}[h!]
\centering
\renewcommand{\arraystretch}{1} 
\begin{tabular}{|l|c|c|} \hline
\textbf{Method} & \textbf{PSNR} & \textbf{SSIM} \\ \hline
NN tanh ($n=15$) & 23.6397 & \textbf{0.4293} \\ \hline
Bilinear Interpolation & \textbf{24.9382} & 0.3921 \\ \hline
Bicubic Interpolation & 24.9244 & 0.3882 \\ \hline
\end{tabular}
\caption{Comparison of rescaling methods for the image \textit{Berlin} (Figure \ref{section_reference_images} (b)).}
\label{Berlintab}
\end{table}
\vspace{0.3cm}
\begin{table}[h!]
\centering
\renewcommand{\arraystretch}{1} 
\begin{tabular}{|l|c|c|} \hline
\textbf{Method} & \textbf{PSNR} & \textbf{SSIM} \\ \hline
NN tanh ($n=15$) & 25.8212 & \textbf{0.3544} \\ \hline
Bilinear Interpolation & \textbf{27.0345} & 0.3245 \\ \hline
Bicubic Interpolation & 27.0127 & 0.3201 \\ \hline
\end{tabular}
\caption{Comparison of rescaling methods for the image \textit{Lisbon} (Figure \ref{section_reference_images} (c)).}
\label{Lisbontab}
\end{table}
\vspace{0.3cm}
\begin{table}[h!]
\centering
\renewcommand{\arraystretch}{1} 
\begin{tabular}{|l|c|c|} \hline
\textbf{Method} & \textbf{PSNR} & \textbf{SSIM} \\ \hline
NN tanh ($n=15$) & 21.3655 & \textbf{0.3516} \\ \hline
Bilinear Interpolation & \textbf{22.8296} & 0.3063 \\ \hline
Bicubic Interpolation & 22.8194 & 0.3024 \\ \hline
\end{tabular}
\caption{Comparison of rescaling methods for the image \textit{Granada} (Figure \ref{section_reference_images} (d)).}
\label{Granadatab}
\end{table}
\FloatBarrier
As evident from the numerical data, Alg. \ref{alg2} turns to be highly effective, particularly in terms of SSIM, when compared to other classical image rescaling methods such as bilinear and bicubic interpolation. This observation holds true across all the selected images. While bilinear interpolation outperforms the NN algorithm in terms of PSNR (see Figure \ref{istogrammi}), the Alg. \ref{alg2} remains a competitive and effective approach overall.

\section{Future works and conclusions}
In this work, we recall some theoretical results on NN operators in the multidimensional framework, including convergence and some quantitative estimates. Two algorithms were implemented: the first focused on modeling remote sensing data, and the second aimed at rescaling or enhancing these data. Experimental results, obtained by applying the algorithms (Alg. \ref{alg1} and Alg. \ref{alg2}) to satellite images from the RETINA dataset, showed that the NN-based methods are particularly effective in terms of SSIM compared to other classical interpolation methods. On the other hand, as observed in Section \ref{CC}, the computational complexity is relatively high and should be optimized. Concerning the effectiveness of the proposed algorithms, we can observe that the good performances in term of SSIM can be justified as in \cite{Marchetti2022}, where it is theorically showed that the converge of the SSIM to $1$ (actually the convergence to zero of the so-called dissimilarity index defined by a continuous version of the SSIM) can be estimated by a quantity depending by the following well-known constant:
\be \label{cf}
c_f\ :=\ {4 \over \sigma_f^2 + c_2}\ +\ {1 \over \mu_f^2 + c_1},
\ee
where $c_1$ and $c_2$ are suitable positive stabilization constants chosen as in \cite{Occorsio2023}, while $\sigma_f$ and $\mu_f$ are respectively the variance and the mean of the image $f$. In particular, from Theorem 3.3 of \cite{Marchetti2022}, when the image $f$ is modeled by ${\cal A}$ defined in (\ref{calA}), we can state the following upper bound
\begin{equation}\label{stima}
|1-\text{SSIM}(K^d_n({\cal A}, \cdot),{\cal A})|\leq c_{\cal A} \left\|K^d_n({\cal A}, \cdot)-{\cal A}(\cdot)\right\|_2^2,
\end{equation}
with ${\cal A} \in L^2(I^d)$, where $SSIM$ here detotes its continuous version. Thus, applying Theorem \ref{theorem3} with $p=2$, we deduce that the right-hand side of (\ref{stima}) converges to zero as $n \to +\infty$.
\\
Note that, the values of $c_f$ in (\ref{cf}) decreases as both $\sigma_f$ and $\mu_f$ increase, and moreover, the bigger contribution to $c_f$ is given by the term with the variance. Indeed, in the images considered for the numerical experiments of this paper, the better results in term of SSIM performances have been reached for {\em Berlin}, which is the image with bigger variance and lower $c_f$ (see Tab. \ref{tab_valori_var}) while the worst accuracy is reached by {\em Granada}, which has the lower variance and the bigger $c_f$, according to what discussed above.
\begin{table}[h!]
\centering
\renewcommand{\arraystretch}{1} 
\begin{tabular}{|l|c|c|c|} \hline
\textbf{Images} & \textbf{Variance} & \textbf{Mean} & \textbf{$c_f$} \\ \hline
\textit{Berlin}  & $9.20\times10^{7}$ & $160.2911$ & $3.89\times10^{-5}$ \\ \hline
\textit{Lisbon}  & $387.6393$        & $195.8720$ & $5.27\times10^{-5}$ \\ \hline
\textit{Rome}    & $1.21\times10^{7}$  & $53.5784$  & $3.47\times10^{-4}$ \\ \hline
\textit{Granada} & $570.6721$          & $76.0948$  & $1.85\times10^{-4}$ \\ \hline
\end{tabular}
\caption{Variance, mean and constants $c_f$ related to the considered images, ordered from the image with the largest $c_f$ to the smallest. The constants $c_f$ are computed using $c_1=0.03\times255$ and $c_2=0.01\times 255$.}
\label{tab_valori_var}
\end{table}

Future work will integrate the multidimensional NN operators approach with Bayesian inversion methods using advanced Monte Carlo Markov Chain (MCMC) techniques. The latter, including Metropolis, Metropolis-Hastings, and Gibbs sampling, allow efficient sampling from the posterior distribution through single flip dynamics \cite{hastings70}. This integration aims at enhancing the accuracy and robustness of the inversion process, improving the retrieval of geophysical variables.
\\
A further study is in preparation to deepen the theoretical understanding of the NN algorithms, with more advanced results on their convergence and approximation order.
\\
As a further interesting open problem, we can mention the possible study of quantitative estimates of the approximation error in terms of a continuous version of the PSNR, in analogy with the theoretical result for the continuous SSIM recalled in (\ref{stima}).

%

\section*{Acknowledgments}

{\small The authors thank the Referees for their constructive comments, which have helped to improve the quality and clarity of the manuscript.
\\
The author are members of the Gruppo Nazionale per l'Analisi Matematica, la Probabilit\`a e le loro Applicazioni (GNAMPA) of the Istituto Nazionale di Alta Matematica (INdAM), of the network RITA (Research ITalian network on Approximation), and of the UMI (Unione Matematica Italiana) group T.A.A. (Teoria dell'Approssimazione e Applicazioni). 
}

\section*{Funding}

{\small The authors have been supported within the project PRIN 2022 PNRR: ``RETINA: REmote sensing daTa INversion with multivariate functional modeling for essential climAte variables characterization", funded by the European Union under the Italian National Recovery and Resilience Plan (NRRP) of NextGenerationEU, under the Italian Ministry of University and Research (Project Code: P20229SH29, CUP: J53D23015950001).

}
\section*{Conflict of interest/Competing interests}

{\small The authors declare that he has no conflict of interest and competing interest.}

\section*{Availability of data and material and Code availability}

{ \small Not applicable.}

\section*{Copyright}
{\small The images (Rome, Berlin, Lisbon, Granada) are contained in the RETINA dataset in \url{https://retina.sites.dmi.unipg.it/dataset.html}. Permission to use, copy, or modify this dataset and its documentation for educational and research purposes only and without fee is granted, provided that this copyright notice and the original authors' names appear on all copies and supporting documentation. This dataset shall not be modified without first obtaining the permission of the authors. The authors make no representations about the suitability of this dataset for any purpose. It is provided "as is" without express or implied warranty.
\\

\noindent In case of publishing results obtained utilizing this dataset, please refer to the following website:
\\
\noindent \url{https://retina.sites.dmi.unipg.it/dataset.html}.



\end{document}